\documentclass[11pt]{article}

\usepackage{latexsym}
\usepackage{amsfonts}
\usepackage[all]{xy}
\usepackage[notref,notcite]{showkeys}

\newtheorem{theorem}{\bf Theorem}[section]

\newtheorem{corollary}[theorem]{\bf Corollary}

\begin{document}

\title{Estimates of the first eigenvalue\\  of minimal hypersurfaces of $\mathbb{S}^{n+1}$}

\author{Abd\^{e}nago Barros \and G. Pacelli Bessa\thanks{\small bessa@mat.ufc.br,
abbarros@mat.ufc.br }}
\date{\today}
\maketitle
\begin{abstract}We consider a solution $f$ of a certain Dirichlet
problem on a domain in $\mathbb{S}^{n+1}$ whose boundary is a
minimal hypersurface and we prove a Poincar\'{e} type inequality
for $f$. Moreover we have an estimate for the first nonzero
eigenvalue for the closed eigenvalue problem on the boundary.
\end{abstract}\section{Introduction}In this note we will let
$M^{n}$ be an embedded compact orientable minimal hypersurface of
$\mathbb{S}^{n+1}$. Yau conjectures that the first nonzero
eigenvalue $\lambda_{1}(M)$ for the closed eigenvalue problem
$\triangle_{M}u +\lambda \,u=0$ on $M$ was equal to $n$. Observe
that $M$ divided $\mathbb{S}^{n+1}$ into two connected components
$\Omega_{1}$ and $\Omega_{2}$ such that $\partial
\Omega_{1}=\partial \Omega_{2}=M$. Choi-Wang \cite{kn:CW} with a
clever idea, applied Reilly  formula to the solution of the
following Dirichlet problem,\begin{equation}\label{eq1}\left\{
\begin{array}{lcl} \overline{\triangle}f&=&0
\,\,\,on\,\,\Omega_{1}\\
& & \\
f&=& \varphi \,\,\,on\,\,M,
\end{array}\right.
\end{equation}where $\varphi$ is the first eigenfunction for the
closed eigenvalue problem on $M$ to prove that
$\lambda_{1}(M)>n/2$. We improve (conceptually) Choi-Wang's
estimates in terms of the soluition $f$ of the problem
(\ref{eq1}), (see Corollary \ref{corollary1}, inequality
\ref{eq3}) with possibility to set up Yau's conjecture provided
one proves equality in \ref{eq5}. The symbols
$\overline{\triangle}$ and $\overline{\nabla}$ will be
respectively  the Laplacian and gradient of the metric of
$\mathbb{S}^{n+1}$ on $\Omega_{1}$ while $\triangle$ and $\nabla$
will be the Laplacian and gradient oh the induced metric on $M$.
\begin{theorem}\label{thm1}Let $f$ be the solution of the Dirichlet problem
(\ref{eq1}). Then \begin{equation}\begin{array}{lcl}p(t)&=&
(2\lambda_{1}(M)-n)\Vert \overline{\nabla} f\Vert^{2}\cdot t^{2}\\
& & \\
& +& 2 \lambda_{1}(M)\Vert f\Vert^{2}\cdot t +
\displaystyle\frac{n}{n+1}\Vert f\Vert^{2} \geq 0,\,\,\,\,\,\,\,
\forall t\in \mathbb{R}
\end{array}
\end{equation}where $\Vert \,\Vert$ denotes the $L^{2}$ norm
on $\Omega_{1}$.
\end{theorem}
\begin{corollary}\label{corollary1}Let $M$ be an orientable
embedded minimal hypersurface of $\mathbb{S}^{n+1}$ and
$\lambda_{1}(M)$ its first non-zero eigenvalue of the Laplacian
for closed eigenvalue problem on $M$. Consider the problem
(\ref{eq1}) and $f$ its solution. Then
\begin{equation}\label{eq3} \lambda_{1}(M) \geq
\displaystyle\frac{n}{2}+\frac{n}{2}\rho(f)
\end{equation}where \begin{equation}\label{eq4}\rho (f)=\displaystyle\frac{2\Vert
\overline{\nabla} f\Vert^{2}-(n+1)\Vert f\Vert^{2}-  2 \Vert
\overline{\nabla} f\Vert^{2}\sqrt{1-(n+1)\frac{\Vert
f\Vert^{2}}{\Vert \overline{\nabla} f\Vert^{2}}}}{(n+1)\Vert
f\Vert^{2}}\end{equation}
\end{corollary}
Observe that $0<\rho (f) \leq 1$ and $\rho(f)=1$ iff $\Vert
\overline{\nabla} f\Vert^{2}=(n+1) \Vert f\Vert^{2}$. Although the
function $f$ does not belong to $H_{0}^{1}(\Omega_{1})$ we have
the following Poincar\'{e} type inequalities.
\begin{corollary}Let $f$ be the solution of the Dirichlet problem
(\ref{eq1}). Then $f$ satisfies the following inequalities:
\begin{equation}\Vert \overline{\nabla}f\Vert^{2}\geq (n+1)\Vert f
\Vert^{2}\label{eq5}
\end{equation}
\begin{equation} \Vert \overline{D}^{\,2} f\Vert^{2}>
\displaystyle\frac{n(n+1)}{4}\Vert f \Vert^{2},
\end{equation}where $ \overline{D}^{\,2} f$ is the Hessian of $f$.
\end{corollary}
\section{Prof of the Results}
Let $\Omega $ be a Riemannian manifold of dimension $n$ with
smooth boundary $\partial \Omega$ and let $f$ be a function on
$\Omega$ which is smooth up to the boundary $\partial \Omega$. We
let $\varphi = f\mid \partial \Omega$ and $u= \partial f/\partial
\nu$ the normal outward derivative of $f$. For $X,Y\in T\Omega$,
$(\overline{D}^{\,2}f)(X,Y)$ denotes  the Hessian tensor. Let
$B(v,w)$ be the second fundamental form of $\partial \Omega$
relative to $\Omega$. Here $v,w$ are tangent to $\partial \Omega$,
$H$ is the mean curvature of $\partial \Omega$ and $ric$ is the
Ricci curvature of $\Omega$. The following identity is known as
the Reilly formula.
\begin{equation}\begin{array}{lcl}\int_{\Omega}
(\overline{\triangle} f)^{2} & = & \int_{\Omega} \vert
\overline{D}^{\,2}f \vert + \int_{\Omega} Ric (\overline{\nabla}f,
\overline{\nabla}f)+ \int_{\partial \Omega}2u\triangle \varphi\\
&& \\
& & +\int_{\partial \Omega}B(\nabla \varphi, \nabla \varphi
)+\int_{\partial \Omega}n H u^{2}\end{array}
\end{equation}
Now we can show the proof of Theorem (\ref{thm1}). If $t=0$ we are
done. Now, for $t\neq 0$ we consider the following Dirichlet
problem\begin{equation}\left\{
\begin{array}{lcl} \overline{\triangle}g&=&f
\,\,\,on\,\,\Omega_{1}\\
& & \\
g&=& t\varphi \,\,\,on\,\,M,
\end{array}\right.
\end{equation}Applying the Green formula we obtain
\begin{equation}\label{eq9}\left\{
\begin{array}{lcl}\int_{M}\varphi \displaystyle\frac{\partial f}{\partial
\nu}&=& \int_{\Omega_{1}}\vert \overline{\nabla}f\vert^{2}\\
&& \\
t\int_{M}\varphi \displaystyle\frac{\partial f}{\partial \nu}&=&
\int_{\Omega_{1}}\langle \overline{\nabla}f, \overline{\nabla}g
\rangle\\
&& \\
\int_{M}\varphi \displaystyle\frac{\partial g}{\partial
\nu}&=&\int_{\Omega_{1}}f^{2} + \int_{\Omega_{1}}\langle
\overline{\nabla}f, \overline{\nabla}g \rangle
\end{array}\right.
\end{equation}
From(\ref{eq9}) we get
\begin{equation}\label{eq10}
\int_{\Omega_{1}}\langle \overline{\nabla}f, \overline{\nabla}g
\rangle=t\int_{\Omega_{1}}\vert
\overline{\nabla}f\vert^{2}\end{equation}and by Cauchy-Schwarz
inequality we get
\begin{equation}\label{eq11}\int_{\Omega_{1}}\vert
\overline{\nabla}g\vert^{2}\geq t^{2}\int_{\Omega_{1}}\vert
\overline{\nabla}f\vert^{2}
\end{equation}from the third equation in (\ref{eq9}) and
(\ref{eq10}) we have that \begin{equation}\label{eq12}
t\int_{M}\varphi \displaystyle\frac{\partial g}{\partial
\nu}=t\int_{\Omega_{1}}f^{2} + t^{2}\int_{\Omega_{1}}\vert
\overline{\nabla}f\vert^{2}.\end{equation}Applying Reilly formula
to $g$, using the fact that $\displaystyle \vert
\overline{D}^{\,2} g\vert^{2}\geq
\frac{1}{n+1}(\overline{\triangle}g)^{2}$ and the assumption that
$\smallint_{M}B(\nabla \varphi, \nabla \varphi )\geq 0$ we have
\begin{equation}\label{eq13}\frac{n}{n+1}\int_{\Omega_{1}}(\overline{\triangle}g)^{2}
\geq n\int_{\Omega_{1}} \vert \overline{\nabla}g \vert^{2}+ 2
\int_{M}\frac{\partial g}{\partial \nu}(\triangle t\varphi)
\end{equation}On the other hand, taking in account  (\ref{eq11}), (\ref{eq12}) and
$\overline{\triangle}g=f$ we have that (\ref{eq13}) implies that
\begin{equation}\frac{n}{n+1}\int_{\Omega_{1}}f^{2}
\geq nt^{2}\int_{\Omega_{1}} \vert \overline{\nabla}f \vert^{2}-
2\lambda_{1}(M)\left[ t\int_{\Omega_{1}}f^{2}
+t^{2}\int_{\Omega_{1}} \vert \overline{\nabla}f \vert^{2}\right]
\end{equation}Therefore we have
\begin{equation}p(t)=(2\lambda_{1}(M)-n)\Vert \overline{\nabla}f
\Vert^{2}\,t^{2} + 2 \lambda_{1}(M) \Vert f \Vert^{2}\, t +
\frac{n}{n+1}\Vert f \Vert^{2}\geq 0.
\end{equation}This finishes the proof of Theorem (\ref{thm1}). The
discriminant of $p$ is non-positive. This can be read as follows
\begin{equation}\label{eq16}(2\lambda_{1}(M)-n) \geq
\frac{n+1}{n}\lambda_{1}(M)^{2} \frac{\Vert f \Vert^{2}}{\Vert
\overline{\nabla}f \Vert^{2}}
\end{equation} From (\ref{eq16}) and (\ref{eq3}) we have the
following poincar\'{e} inequality for $f$,
\begin{equation}\Vert
\overline{\nabla}f \Vert^{2}\geq (n+1)\Vert f \Vert^{2}
\end{equation}In the proof of Theorem (\ref{thm1}) we did not
count on with an extra term $ \smallint_{M} B(\nabla t\varphi,
\nabla t\varphi)\geq 0$ on the right side of (\ref{eq13}). Taking
it in account we have in fact that $p(t) \geq \smallint_{M}
B(\nabla \varphi, \nabla \varphi)\cdot t^{2}$. From that we can
conclude that $$(2\lambda_{1}(M)-n)\Vert \overline{\nabla}f
\Vert^{2}- \smallint_{M} B(\nabla \varphi, \nabla \varphi)\geq
\frac{n+1}{n}\lambda_{1}(M)^{2} \frac{\Vert f \Vert^{2}}{\Vert
\overline{\nabla}f \Vert^{2}}$$On the other hand, Reilly formula
also gives $$(2\lambda_{1}(M)-n)\Vert \overline{\nabla}f
\Vert^{2}=\Vert \overline{D}^{2}f \Vert^{2}+\smallint_{M} B(\nabla
\varphi, \nabla \varphi).$$ Therefore we obtain
\begin{equation}\Vert \overline{D}^{2}f \Vert^{2}\geq
\frac{n+1}{n}\lambda_{1}(M)^{2}\Vert
f\Vert^{2}>\frac{n(n+1)}{4}\Vert f\Vert^{2},
\end{equation}since $\lambda_{1}(M)>n/2$.

\end{document}